\documentclass[11pt]{article}

\usepackage{latexsym}
\usepackage{a4}
\usepackage{amsmath}
\usepackage{amsfonts}
\usepackage{amssymb}
\usepackage{amscd}
\usepackage{graphicx}
\usepackage{graphics}
\usepackage{fancyhdr}
\usepackage{makeidx}

\newcommand{\dis}{\displaystyle}
\newtheorem{thm}{Theorem}[section]
\newtheorem{prop}[thm]{Proposition}
\newtheorem{lem}[thm]{Lemma}

\newtheorem{Def}[thm]{Definition}

\newenvironment*{Proof}{{\bf Proof.}}

\newcommand{\bbb}[1]{\mbox{\boldmath$#1$}}

\newcommand{\fa}{\forall}

\newcommand{\ra}{\;\rightarrow\;}

\newcommand{\al}{\alpha}
\newcommand{\bi}{\beta}

\newcommand{\Ga} {{\varGamma}}
\newcommand{\Sig} {{\varSigma}}
\newcommand{\de}{\delta }
\newcommand{\OO} {{\varOmega}}

\newcommand{\e}{\varepsilon }

\newcommand{\Fi}{\varPhi}

\newcommand{\zi}{\zeta }
\newcommand{\thi}{\theta }
\newcommand{\vthi}{\vartheta }

\newcommand{\La} {{\mit\Lambda}}

\newcommand{\ti}{\tau }

\newcommand{\C}{\mathbb{C}}
\newcommand{\R}{\mathbb{R}}

\newcommand{\ssum}{\sum\limits}

\newcommand{\oD}{\overline{D}}
\newcommand{\oO}{\overline{\varOmega}}

\newcommand{\td}{\widetilde{d}}

\newcommand{\tV}{\widetilde{V}}
\newcommand{\tA}{\widetilde{A}}
\newcommand{\tM}{\widetilde{M}}

\newcommand{\cD}{{\cal{D}}}

\newcommand{\ld}{\ldots}

\newcommand{\sm}{\smallsetminus}

\newcommand{\qb}{$\quad\blacksquare$}

\begin{document}
\title{\bf An extension of the disc Algebra}
\author{V. Nestoridis}
\date{}
\maketitle
\begin{abstract}
We identify all uniform limits of polynomials on the closed unit
disc $\overline{D}$ with respect to the chordal metric $\chi$ on
$\C\cup\{\infty\}$. One such limit is $f\equiv\infty$. The other
limits are holomorphic functions $f:D\ra\C$ so that for every
$\zi\in\partial D$ the $\dis\lim_{z\ra \zi\atop z\in D}f(z)$
exists in $\C\cup\{\infty\}$. The class of the above functions is
denoted by $\tA(D)$. We study properties of the members of
$\tA(D)$, as well as, some topological properties of $\tA(D)$
endowed with its natural metric topology. There are several open
questions and new directions of investigation.
\end{abstract}
{\em AMS classification number}: Primary 30J99, secondary 46A99,
30E10 \vspace*{0.1cm} 

\section{Introduction}  
\noindent

In this article we consider $\chi$ the chordal metric in
$\C\cup\{\infty\}$, where $\C$ denotes the complex plane. First we
investigate the set of uniform limits on the closed unit disc
$\oD$ of the polynomials with respect to the metric $\chi$. One
such limit is $f\equiv\infty$. The other limits are exactly
functions $f:D\ra C$ holomorphic in the open unit disc $D$, such
that for every $\zeta\in T=\partial D$ the limit $\dis\lim_{z\ra
\zi\atop z\in D}f(z)$ exists in $\C\cup\{\infty\}$. The class of
the above functions is an extension of the disc algebra $A(D)$ and
is denoted by $\tA(D)$.

If $f\in \tA(D)$ is such that $f\not\equiv\infty$, then Privalov's
Theorem implies that $f^{-1}(\infty)$ is a compact subset of the
unit circle $T$ with zero length. Reversely for every such set,
there exist $f\in\tA(D)$, $f\not\equiv\infty$ so that
$f^{-1}(\infty)$ is exactly this set. Compact subsets of $T$ with
positive length are not compact of interpolation for $\tA(D)$. It
is an open question if every compact subset $E\subset T$ with zero
length is a compact of interpolation for $\tA(D)$. That is, is it
true that every continuous function $h:E\ra\C\cup\{\infty\}$ has
an extension in $\tA(D)$. For general $f,g\in\tA(D)$ it is not
true that
$\dis\sup_{z\in\oD}\chi(f(z),g(z))\le\dis\sup_{|z|=1}\chi(f(z),g(z))$,
even if we multiply the right member by a constant $C<+\infty$
independent of $f,g$. However two members of $\tA(D)$ which
concide on $T$, should be identical. The mean value property fails
in general in $\tA(D)$.

The disc algebra $A(D)$ endowed with the supremum norm is a Banach
algebra, so in particular a complete metric space.  The larger
class $\tA(D)$ is naturally endowed with the metric
$\td(f,g)=\dis\max_{|z|\le1}\chi(f(z),g(z))$, $f,g\in\tA(D)$. The
metric space $(\tA,\td)$ is also complete. Often we write
$\chi(f,g)$ for $\td(f,g)$, $f,g\in\tA(D)$. Then $A(D)$ is an open
dense subset of $\tA(D)$. Furthermore the relative topology of
$A(D)$ from $\tA(D)$ coincides with the usual topology of $A(D)$.
For $f\in\tA(D)$ we denote by $E_f$ the set $E_f=\{\zi\in
T:f(\zi)\notin f(D)\}$. We prove that generically, for every
$f\in\tA(D)$ the set $E_f$ has zero length. Also we show that the
set $Y=\{f\in\tA(D):f(D)\subset f(T)\}$ is a non-void closed
subset of $\tA(D)$ of first category. We also introduce the set
$W=\{f\in\tA(D):f(T)=\C\cup\{\infty\}\}$. This is a closed subset
of $\tA(D)$ of first category, but we do not know if it is
non-empty. If we assume that every compact subset of $T$ with zero
length is a compact of interpolation for $\tA(D)$, then we can
show that $W$ is non-empty.

Another open question is to characterize the compact subsets of
$\oD$ which are of interpolation for $\tA(D)$. A necessary condition
is that, they are of the form $E\cup\{z_n:n\in I\}$, $I$ finite or
infinite denumerable, $|z_n|<1$ for all $n\in I$, where $E\subset T$
is compact with zero length and all limit points of $z_n$, $n\in I$
are included in $E$. Is the same condition sufficient or not? If
not, what is a nessary and sufficient condition? Another open
question is to characterize the zero sets in $\oD$ of the members of
$\tA(D)$. Every result which holds in $A(D)$ may be examined if it
holds in $\tA(D)$.

In the above discussion the uniform convergence with respect to
the metric $\chi$ was considered mainly on $\oD$, but we can
consider it on other compact subsets of $\C$. For instance, does
Mergelyan's Theorem holds for this uniform convergence with
respect to the metric $\chi$? In general, this is an open
question:\vspace*{0.2cm} \\
{\bf Question:} Let $L\subset\C$ be a compact set with connected
complement. Let $f:L\ra\C\cup\{\infty\}$ be a continuous function
such that for every component $V$ of $L^0$ either
$f_{|V}\equiv\infty$ or $f(V)\subset\C$ and $f_{|V}$ is
holomorphic. Is then true that there exists a sequence of
polynomials converging to $f$ uniformly on $L$ with respect to the
metric $\chi$?

It is easy to see that the converse holds.

In the present article we give a positive answer to the previous
question in a few particular cases, while the general case is
open.

More generally one can investigate the uniform limits with respect
to the metric $\chi$ on other sets, as an annulus or a circle, of
a sequence of rational functions with prescribed set of poles. For
instance I do not know what are the limits on a circle when we
approximate by polynomials. On the contrary on the unit circle, if
we approximate with trigonometric polynomials, the limits are all
continuous functions $f:T\ra\C\cup\{\infty\}$.

Therefore, we see that there are several directions to continue
our investigation.

Before closing we say that, if we consider $H(D)$, the space of
holomorphic functions in $D$, endowed with the topology of uniform
convergence on compacta, then $\tA(D)\cap H(D)$ is a subset of
$H(D)$ of first category. One reason is that $\tA(D)$ is disjoint
with the residual set $U$ in $H(D)$ of universal Taylor series
(\cite{1}, \cite{2}). In fact if $f\in U$, then it is impossible
to have a finite (in $\C$) limit $\dis\lim_{z\ra
e^{i\vthi}\atop|z|<1}f(z)$ (\cite{3}).
\section{The definition} 
\noindent

We consider the open unit disc $D=\{z\in\C:|z|<1\}$ in the complex
plane $\C$ and we fixe $w\in\C\cup\{\infty\}$. We are looking for
the usual uniform limits on the closed unit disc
$\oD=\{z\in\C:|z|\le\}$ of rational functions whose poles in
$\C\cup\{\infty\}$ are included in $\{w\}$. For $w=\infty$ we are
looking for the usual uniform limits on $\oD$ of polynomials
$P(z)=\sum\limits^N_{n=0}a_nz^n$, $a_n\in\C$, $N=0,1,2,\ld\,.$ It is
well known that the set of such limits coincides with the disc
Algebra $A(D)=\{f:\oD\ra\C$, continuous on $\oD$ and holomorphic in
$D\}$. If $w\in\C$, $|w|>1$ we are looking for the usual uniform
limits on $\oD$ of functions of the form
$P(z)=\sum\limits^N_{n=0}a_n\frac{1}{(z-w)^n}$, $a_n\in\C$,
$N=0,1,2,\ld\,.$ Then the set of such limits ia again the disc
Algebra $A(D)$.

If $w\in\C$, $|w|\le1$, the function $P(z)=f(z)$ satisfies
$f(w)=\infty$, with $w\in\oD$ or $f$ is constant. Looking for
uniform limits with respect to the Eucledian metric in
$\C\cong\R^2$ we exclude the cases $f(w)=\infty$, so the set of
such limits should be the set of constant functions $g(z)=c$ for
all $z\in\oD$ with $c\in\C$. The answer is not satisfactory,
because the set of limits does not contain the trivial cases of
constant sequences as $f_n(z)=\frac{1}{z-w}$, $n=1,2,\ld\,.$ Such
a sequence should converge to the function $f(z)=\frac{1}{z-w}$;
but, this is excluded, because it takes the value $\infty$ in
$z=w\in\oD$. This leads us to search for the uniform limits of the
rational functions with poles in $\C\cup\{\infty\}$ included in
$\{w\}$ with respect to the chordal distance $\chi$ on
$\C\cup\{\infty\}$.

We have
\[
\chi(z_1,z_2)=\frac{|z_1-z_2|}{\sqrt{1+|z_1|^2}\sqrt{1+|z_2|^2}} \
\ \text{for} \ \ z_1,z_2\in\C
\]
\[
\text{and} \ \ \chi(z,\infty)=\frac{1}{\sqrt{1+|z|^2}} \ \
\text{for} \ \ z\in\C \ \ \text{and} \ \ \chi(\infty,\infty)=0.
\]
Geometrically, we identify $\C\cup\{\infty\}$ with
$S^2\subset\R^3$ via stereographic projection and $\chi$ is a
constant multiple of the restriction on $S^2$ of the usual
Eucledian metric in $\R^3$.

If $f_n,f:E\ra\C\cup\{\infty\}$, $n=1,2,\ld$ are functions defined
on a set $E$, then $f_n\ra f$ uniformly on $E$, with respect to
the metric $\chi$ if and only if, $\dis\sup_{z\in
E}\chi(f_n(z),f(z))\ra0$, as $n\ra+\infty$. With this convergence
in mind, we must reexamine the cases $w=\infty$ or $w\in\C$,
$|w|>1$, as well.

Suppose $w=\infty$. Assume that a sequence of polynomials $f_n$
$n=1,2,\ld$ converges to a function $f:\oD\ra\C\cup\{\infty\}$
uniformly on $\oD$ with respect to $\chi$. Since uniform
convergence preserves continuity, it follows that $f$ is
continuous when $\C\cup\{\infty\}$ is endowed with the metric
$\chi$. Suppose that $f(z_0)=\infty$ for some $z_0$, $|z_0|<1$.
Then for some $r>0$ small we have
$\chi(f(z_0)=\infty,f(z))<\frac{1}{3}$ for all $z:|z-z_0|<r$.
Since $f_n\ra f$ uniformly with respect to $\chi$, there exists
$n_0$ so that $\chi(f_n(z),\infty)<\frac{1}{2}$ for all $n\ge n_0$
and all $z:|z-z_0|<r$. Thus, $f_n(z)\in\C-\{0\}$ for all $n\ge
n_0$ and $z:|z-z_0|<r$. Since $f_n$ is holomorphic in $\C$, being
a polynomial, it follows that $\frac{1}{f_n}$ is holomorphic in
the disc $\{z:|z-z_0|<r\}$. As
$\chi\Big(\frac{1}{z_1},\frac{1}{z_2}\Big)=\chi(z_1,z_2)$ for all
$z_1,z_2\in\C\cup\{\infty\}$, it follows that
$\frac{1}{f_n}\ra\frac{1}{f}$ uniformly on $\{z:|z-z_0|<r\}$ with
respect to the distance $\chi$. Since
$\chi(\infty,f(z))<\frac{1}{3}$ it follows
$\chi\Big(0,\frac{1}{f(z)}\Big)<\frac{1}{3}$ for all $z:|z-z_0|<r$
and $\chi\Big(0,\frac{1}{f_n(z)}\Big)<\frac{1}{2}$ for all
$z:|z-z_0|<4$ and $n\ge n_1$, for some $n_1\ge n_0$. Thus, there
exist $M<+\infty$ so that $\Big|\frac{1}{f_n(z)}\Big|<M$ and
$\Big|\frac{1}{f(z)}\Big|<M$ for all $n\ge n_1$ and $z:|z-z_0|<r$.
It follows easily
\[
\sup_{z:|z-z_0|<r}\bigg|\frac{1}{f_n(z)}-\frac{1}{f(z)}\bigg|\le\sup_{z:
|z-z_0|<r}\chi\bigg(\frac{1}{f_n(z)},\frac{1}{f(z)}\bigg)\cdot C \
\ \text{for all} \ \ n\ge n_1 \ \ \text{with} \ \ C<+\infty.
\]
Thus, the sequence of holomorphic mappings $\frac{1}{f_n}$,
$n\ge1$ converges uniformly on $\{z:|z-z_0|<r\}$ to $\frac{1}{f}$
with respect to the usual Eucledian distance on $\C=\R^2$. But
$\frac{1}{f_n(z)}\neq0$ for all $z:|z-z_0|<r$, $n\ge n_1$ and
$\frac{1}{f(z_0)}=0$. Hurwitz Theorem implies $f\equiv\infty$ on
$\{z:|z-z_0|<r\}$. Thus the set $f^{-1}(\infty)\cap D$ is open.
The same set is relatively closed by continuity. Since $D$ is
connected, it follows $f\equiv\infty$. Reversely, for
$P_n(z)\equiv n$ the limit is $\infty$. So the constant infinity
function is one possible limit. If $f$ is not identically equal to
$\infty$ and $f$ is a possible limit, then $f(D)\subset\C$ and the
value $\infty$ is possible only on the boundary $\partial D=T$
where $T=\{z\in\C:|z|=1\}$ is the unit circle. It is also easily
seen, by the previous discussion, that $f_n\ra f$ uniformly on
compact subsets of $D$ with respect to the usual Eucledian metric
in $\C=\R^2$. Thus, $f$ is holomorphic in $D$.
\begin{thm}\label{thm2.1} 
Let $f:D\ra\C\cup\{\infty\}$ be a continuous function, such that,
$f(D)\subset\C$ and $f_{|D}:D\ra C$ is a holomorphic function.
Then, there exists a sequence of polynomials
$f_n(z)=\sum\limits^{N_n}_{j=0}a^n_jz^j$, $a^n_j\in\C$,
$N_n\in\{0,1,2,\ld\}$ so that $f_n\ra f$ uniformly on $\oD$ with
respect to the metric $\chi$.
\end{thm}
\begin{Proof}
Let $\e>0$. Continuity of $f$ on the compact set $\oD$ implies
uniform continuity with respect to the metric $\chi$. Thus, there
exists $r$, $0<r<1$, so that $\chi(f(z),f(rz))<\frac{\e}{2}$ for all
$z\in\oD$. We consider the Taylor development of $f$ in $D$ with
center zero $f(z)=\ssum^\infty_{j=0}a_jz^j$. The convergence is
uniform on $\{z:|z|\le r\}$ with respect to the usual Eucledian
metric in $\C=\R^2$. Thus, there exists $N$ so that
\[
\bigg|f(\ti)-\sum^N_{j=0}a_j\ti^j\bigg|<\frac{\e}{2}, \ \ \text{for all}
 \ \ \ti \ \ \text{with} \ \ |\ti|\le r.
\]
But $\chi(z_1,z_2)\le|z_1-z_2|$ and $|rz|\le r$ for all $z\in\oD$.
It follows
\[
\chi\bigg(f(rz),\sum^N_{j=0}a_jr^jz^j\bigg)<\frac{\e}{2} \ \
\text{for all} \ \ z\in\oD.
\]
The triangular inequality implies $\chi(f(z),P(z))<\e$ for all
$z\in\oD$ where $P$ is the polynomial
$P(z)=\ssum^N_{j=0}a_jr^jz^j$. This completes the proof. \qb
\end{Proof}
\begin{Def}\label{Def2.2}
$\tA(D)$ denotes the set of all functions
$f:\oD\ra\C\cup\{\infty\}$ such that $f$ is identically equal to
$\infty$ of $f$ is continuous with respect to the metric $\chi$,
$f(D)\subset\C$ and $f_{|D}:D\ra\C$ is holomorphic in $D$.
\end{Def}

The class $\tA(D)$ contains the set of polynomials and we saw
previously that polynomials are dense in $\tA(D)$, for the
topology of uniform convergence on $\oD$ with respect to the
metric $\chi$.

Continuing our effort to identify all uniform limits on $\oD$ with
respect to the metric $\chi$ by functions of the form
$P(z)=\ssum^N_{j=0}a_j\frac{1}{(z-w)^j}$, where $w\in\C$ is fixed,
we can say the following.

If $|w|>1$ then the set of limits is $\tA(D)$, the same as for
$w=\infty$. The only modification is in the proof of Theorem
\ref{thm2.1} Instead of approximating $f$ on $\{z:|z|\le r\}$ by a
partial sum of the Taylor development, we can use Runge's Theorem
and make approximation by a polynomial in $\frac{1}{z-rw}$. If
$|w|<1$, then the set of possible limits is the set of
$f\equiv\infty$ or $f(z)=Q\Big(\frac{1}{z-w}\Big)$, where $Q$ is any
polynomial; that is, the set of limits is almost the same as the set
of approximating functions. To prove this, assume
$f\not\equiv\infty$. Then, Hurwitz Theorem implies that
$f(D-\{w\})\subset\C$. We already used this fact which states that
the limiting function $f$ can not take the value $\infty$ at an
interior point where the approximating functions take finite values,
unless $f\equiv\infty$. Since $f$ is holomorphic in $D-\{w\}$, then
there exist unic functions $f_1$ and $f_2$ so that $f=f_1+f_2$ on
$D-\{w\}$, $f_1$ is holomorphic in $D$ and $f_2$ is holomorphic in
$\C\sm\{w\}$ and satisfies $\dis\lim_{z\ra\infty}f_2(z)=0$. This
follows using the Laurent decomposition of $f$ with center $w$; see
also \cite{4}. If $0<r_1<r_2$ are such that $\{\zi:|\zi-w|\le
r_2\}\subset D$, then $f_1(z)=\frac{1}{2\pi
i}\int\limits_{|\zi-w|=r_2}\frac{f(\zi)}{\zi-z}d\zi$ for $|z-w|\le
r_1$.

Suppose $P_n\Big(\frac{1}{z-w}\Big)\ra f(z)$ uniformly on $\oD$
with respect to $\chi$, where $P_n$ are polynomials. Since $f$
takes finite values on $\{\zi:|\zi-w|=r_2\}$ whose set is a
compact set not containing $\infty$, it follows easily that
$P_n\Big(\frac{1}{\zi-w}\Big)\ra f(\zi)$ uniformly on
$\{\zi:|\zi-w|=r_2\}$ with respect to the usual Euclidean metric
on $\C=\R^2$. It follows that
\[
\frac{1}{2\pi
i}\int_{|\zi-w|=r_2}\frac{P_n\Big(\frac{1}{\zi-w}\Big)}{\zi-z}d\zi\ra
f_1(z) \ \ \text{on} \ \ \{z:|z-w|\le r_1\}.
\]
But the left hand side coincide with the constant term of the
polynomial $P_n$. So $f_1$ is constant on $\{z:|z-w|\le r_1\}$ and
on $\oD$ by analytic continuation and continuity. Thus,
$f_1(z)\equiv c\in\C$. The function $f_2$ coincides with the
principal part of the Laurent development of $f$ with center $w$.
$w$ is a pole or a removable singularity for $f$. If $w$ is a pole
for $f$, we conclude $f_2(z)=\ssum^N_{j=1}b_j\frac{1}{(z-w)^j}$
with $N=\{1,2,\ld\}$, $b_j\in\C$. If $w$ is a removable
singularity for $f$ we conclude $f_2\equiv0$. This completes our
investigation in the case $|w|<1$. It remains the case $|w|=1$. In
this case the functions $f\equiv\infty$ or
$f(z)=P\Big(\frac{1}{z-w}\Big)$ with $P$ any polynomial are among
the limiting functions. If $f$ is a limiting function
$f\not\equiv\infty$, then $f(D)\subset\C$ and the value $\infty$
is
possible only on the unit circle $T$. \vspace*{0.2cm} \\
{\bf Question:} Is there any limiting function $f$ in the case
$|w|=1$, such that $f\not\equiv\infty$ and $f(\zi)=\infty$ for
some $\zi:|\zi|=1$, $\zi\neq w$? What is a characterization of the
set of limiting functions in this case?
\section{Properties of the members of $\bbb{\tA(D)}$}  
\noindent
\begin{prop}\label{prop3.1} 
Let $f\in\tA(D)$ and $c\in\C\cup\{\infty\}$. If $f(\zi)=c$ holds
for a set of positive Lebesgue measure on the unit circle $T$,
then $f\equiv c$.
\end{prop}
\begin{Proof}
If $c\in\C$, then the result follows by Privalov's Theorem
(\cite{5} page 84) applied to the function $f-c$. Let $c=\infty$.
Then, by uniform continuity of $f$ with respect to the metric
$\chi$, there exists $\thi_1<\thi_2<\thi_1+2\pi$ and $r$, $0<r<1$,
so that $f(z)\neq0$ on $\{te^{i\vthi}:r\le
t\le1,\thi_1\le\vthi\le\thi_2\}\equiv L$ and $f(\zi)=\infty$ on a
compact subset $E$ of $\{e^{i\vthi}:\thi_1\le\vthi\le\thi_2\}$
with positive length. Assume $f\not\equiv\infty$ to arrive at a
contradiction.

There exists a Riemann mapping $F:D\ra L^0$ with $F'\in H^1$,
because $\partial L$ is rectifiable (\cite{6} page 44). For the
length of $E$, $|E|>0$, we have
$|E|=\int\limits_{F^{-1}(E)}F'(e^{i\vthi})d\vthi$. Thus,
$F^{-1}(E)$ is a compact subset of $T=\partial D$ with positive
length and the function $\frac{1}{f}\circ F$ vanishes on
$F^{-1}(E)$. Since $\frac{1}{f}\circ F$ is holomorphic in $D$,
Privalov Theorem implies $\frac{1}{f}\circ F\equiv 0$ on $\oD$,
which gives $f\equiv\infty$.  \qb
\end{Proof}

\begin{prop}\label{prop3.2} 
Let $E\subset T=\partial D$ be a compact set with zero length. Then,
there exists $f\in\tA(D)$ so that $f^{-1}(\infty)=E$.
\end{prop}
\begin{Proof}
It is well known (\cite{7} page 81) that there exists $g\in A(D)$,
such that, $g_{|_E}\equiv 1$ and $|g(z)|<1$ for all $z\in\oD-E$.
It suffices to set $f=\frac{1}{g-1}$.   \qb
\end{Proof}

If a compact set $E\subset T=\partial D$ has positive length, then
it is not a compact of interpolation for $\tA(D)$. Let $\zi_0\in
E$. We consider the function $h(\zi_0)=1$, $h(\zi)=0$ for $\zi\in
E\sm\{e^{i\vthi}\zi_0:|\vthi|\le\e\}$, where $\e=\frac{|E|}{100}$,
extended linearly on $\{e^{i\vthi}\zi_0:|\vthi|\le\e\}\cap E$.
This function $h$ can not have an extension $f$ in $\tA(D)$. Since
$|E\sm\{e^{i\thi}\zi_0:|\vthi|\le\e\}|>0$, it would follow that
$f\equiv 0$; but $f(\zi_0)=h(\zi_0)=1\neq0$. This gives a
contradiction. \vspace*{0.2cm} \\
{\bf Question:} Is it true that every compact set $E\subset
T=\vthi D$ with zero length is a compact of interpolation for
$\tA(D)$? That is, is it true, that, for every continuous function
$h:E\ra\C\cup\{\infty\}$ there exists $f\in\tA(D)$ so that
$f_{|E}=h$? Furthermore one could ask for a characterization of
compact sets $E\subset\oD$ which are of interpolation for
$\tA(D)$. \vspace*{0.1cm}

We notice that for $f,g\in\tA(D)$ we do not have in general
$\dis\sup_{z\in\oD}\chi(f(z),g(z))\le\dis\sup_{|z|=1}\chi(f(z),g(z))\cdot
C$ for any constant $C<+\infty$ independent of $f$ and $g$. In
fact, even if we fix $r$, $0<r<1$, there is no constant
$C<+\infty$ so that $\dis\sup_{z\in\oD}\chi(f(z),g(z))\le
C\dis\sup_{r\le|z|\le1}\chi(f(z),g(z))$ for all $f,g\in\tA(D)$. To
see this we set $g\equiv\infty$ and $f_n(z)=nz$.

Then,
\[
\sup_{z\in\oD}\chi(f_n(z),g(z))\ge\chi(0,\infty)=1 \ \ \text{but}
\]
\[
\sup_{r\le|z|\le1}\chi(f_n(z),g(z))=\chi(nr,\infty)\ra0, \ \
\text{as} \ \ n\ra+\infty.
\]
However, we have the following.
\begin{prop}\label{prop3.3} 
Let $f,g\in\tA(D)$. If $f(\zi)=g(\zi)$ for all $\zi\in T=\partial
D$, then $f\equiv g$.
\end{prop}
\begin{Proof}
If $f\equiv\infty$ or $g\equiv\infty$, then the result follows
from Proposition \ref{prop3.1} We assume $f\not\equiv\infty$ and
$g\not\equiv\infty$. Proposition \ref{prop3.1} implies that the
set $E=\{\zi\in T:f(\zi)=\infty\}$ is a compact set with zero
length. Since $f(\zi)=g(\zi)$ for all $\zi\in T$ we conclude that
the function $f-g$, which takes finite values in $D$ and is
holomorphic in $D$, extends continuously on $D\cup(T\sm E)$ with
values in $\C$. We also have $f(z)-g(z)=0$ on $T\sm E$. Since
$T\sm E$ contains a compact set with positive length, it follows
by Privalov's Theorem (\cite{5}, page 84) that $f(z)\equiv g(z)$
on $D$ and by continuity on $\oD$.

The proof is complete. \qb
\end{Proof}

Furthermore the mean value property
$f(0)=\frac{1}{2\pi}\int\limits^{2\pi}_0f(e^{i\vthi}d\vthi$ does
not hold for all $f\in\tA(D)$, $f\not\equiv\infty$. If we set
$\frac{1}{z-1}=f(z)$. Then $f\in\tA(D)$ with
$f^{-1}(\infty)=\{1\}$ and $f\notin L^1(T,d\vthi)$. Thus, the mean
value property is not valid. Even if we interprete the integral as
a principal value we have:
\[
\lim_{\e\ra0^+}\frac{1}{2\pi}\int^{2\pi-\e}_\e\frac{1}{e^{i\vthi}-1}d\vthi=
\lim_{\e\ra0^+}\frac{1}{2\pi}\int^\pi_\e2Re\frac{1}{e^{i\thi}-1}d\vthi=
-\frac{1}{2}\neq-1=f(0),
\]
because $Re\frac{1}{e^{i\vthi}-1}=-\frac{1}{2}$ for all
$\vthi\in\R$.

The mean value property
$f(0)=\frac{1}{\pi}\int\limits_{x^2+y^2<1}f(x+iy)dxdy$ is neither
valid for all $f\in\tA(D)$, $f\not\equiv\infty$. If we set
$f(z)=\frac{1}{(z-1)^2}$ then $f\in\tA(D)$, $f^{-1}(\infty)=\{1\}$
and $f\notin L^1(D,dxdy)$. However, in polar coordinates, the
iterated integral
$\frac{1}{\pi}\int\limits^1_{r=0}\int\limits^{2\pi}_{\thi=0}f(re^{i\vthi})(d\vthi)
\cdot rdr$ is equal to $f(0)$ for all holomorphic functions in
$D$. I also think that for every $G\in H^1_0$ the formula
$f(0)=\frac{1}{2\pi}\int\limits^{2\pi}_{\vthi=0}f(e^{i\vthi})(1+G(e^{i\vthi}))d\vthi$
does not hold for all $f\in\tA(D)$, $f\not\equiv0$, but I do not
have a proof.
\section{Topological properties of $\bbb{\tA(D)}$}
\noindent

In the disc algebra $A(D)=\{f:\oD\ra\C$ continuous on $\oD$ and
holomorphic in $D\}$ we consider the metric
$d(f,g)=\dis\max_{|z|\le1}|f(z)-g(z)|$. Then $(A(D),d)$ is a
complete metric space. In $\tA(D)$ we consider the metric
$\td(f,g)=\dis\max_{|z|\le1}\chi(f(z),g(z))$. Very often, we write
$\chi(f,g)$ instead of $\td(f,g)$, $f,g\in\tA(D)$. Then $(A(D),
\td)$ is a complete metric space. If $f_n\in\tA(D)$ is
$\td$-Cauchy sequence, then, since $(\C\cup\{\infty\},\chi)$ is a
complete metric space, there exists a function $f:\oD\ra
\C\cup\{\infty\}$, such that, for every $z\in\oD$
$\chi(f_n(z),f(z))\ra0$, as $n\ra+\infty$. Let $\e>0$, then
$\dis\sup_{z\in\oD}\chi(f_n(z),f_m(z))<\e$ for all $n,m\ge n_0$,
for some $n_0$. Making $m\ra+\infty$ we obtain $\td(f_n,f)<\e$ for
all $n\ge n_0$. So $\td(f_n,f)\ra0$. If $f\equiv\infty$, then
$f\in\tA(D)$.

Assume $f(z_0)\in\C$ for some $z_0\in\oD$.

Thus $f_n(z_0)\in\C$ for all $n\ge n_1$ for some $n_1$. On a
compact neighborhood $V_{z_0}$ of $z_0$ in $\oD$ we find
$|f(z)|\le M$ for some $M<+\infty$ and for all $z\in\tV_{z_0}$
where $\tV_{z_0}$ is a compact neighborhood of $z_0$ included in
$V_{z_0}$.

We set $\e=\frac{1}{2}\chi(\infty,\{a\in\C:|\al|\le M\})>0$. Then,
the set $\{\bi\in\C\cup\{\infty\}:\chi(\bi,\{a\in\C:|\al|\le
M\})\le\e\}$ is compact and does not contain $\infty$. But every
compact subset of $\C$ is bounded (in Eucledian distance). So
there exists $M'$, $M\le M'<+\infty$, so that,
$\chi(\bi,\al)\le\e$ for some $\al$ with $|\al|\le M$ implies
$|\bi|\le M'$. Thus, there exist $n_1$ so that $|f_n(z)|\le M'$
and $|f(z)|\le M'$ for all $z\in\tV_{z_0}$ and $n\ge n_1$. Thus
$\chi(f_n(z),f(z))\ge\frac{|f_n(z)-f(z)|}{1+M'^2}$. So $f_n\ra f$
uniformly with respect to Eucledian distance on $\tV_{z_0}$. It
follows that $f$ is holomorphic in some open set $U_{z_0}$
containing $z_0$, if $z_0\in D$ or $z_0\in\partial U_{z_0}$ if
$z_0\in T$.

Assume $z_1\in D$ is such that $f(z_1)=\infty$. Then on a compact
neighborhood $U_{z_1}$ of $z_1$ we have
$\chi(f(z),\infty)<\frac{1}{3}$. There exists $n_3$ so that
$\chi(f(z),f_n(z))<\frac{1}{3}$ for all $z\in U_z$, and $n\ge
n_3$. The triangular inequality implies
$\chi(f_n(z),\infty)<\frac{2}{3}$ for all $z\in U_z$, and $n\ge
n_3$. This implies $|f_n(z)|>\sqrt{\frac{7}{2}}$ and
$|f(z)|>\sqrt{\frac{7}{2}}$. Thus, in the interior of $U_{z_1}$
the holomorphic functions $\frac{1}{f_n}$ and $\frac{1}{f}$ are
bounded above by $\sqrt{\frac{2}{7}}$. Since
$\chi\Big(\frac{1}{\al},\frac{1}{\bi}\Big)=\chi(\al,\bi)$ we see
$\frac{1}{f_n}\ra\frac{1}{f}$ uniformly on $U_{z_1}$ with respect
to the metric $\chi$; but
$\chi\Big(\frac{1}{f_n(2)},\frac{1}{f(z)}\Big)\ge\frac{\Big|
\frac{1}{f_n(z)}-\frac{1}{f(z)}\Big|}{\frac{9}{7}}$ on $U_{z_1}$.
This implies that $\frac{1}{f_n}\ra\frac{1}{f}$ uniformly on
$U_{z_1}$ with respect to the usual Eucledian metric on $\C=\R^2$.
If $f_n\equiv\infty$ for infinitely many $n$'s, then obviously
$f\equiv\infty$ belongs to $\tA(D)$. So we assume
$f_n(D)\subset\C$ and $f_{n|D}$ holomorphic for every $n$. But
$\frac{1}{f_n(z)}\neq0$ for every $n$ and $z\in U_{z_1}$ while
$\frac{1}{f(z_1)}=0$. Hurwitz Theorem implies $f(z)\equiv\infty$
on $U_{z_1}$. Thus the set $\{z\in D:f(z)=\infty\}$ is open. By
the continuity of $f$, the same set is closed in $D$. Since $D$ is
connected, this set is either empty or the hole $D$. In the second
case $f\equiv\infty$, so $f\in\tA(D)$. In the first case
$f(D)\subset\C$ and $f$ is holomorphic in $D$. Since
$f:\oD\ra\C\cup\{\infty\}$ is continuous we conclude that
$f\in\tA(D)$. This completes the proof. \qb\vspace*{0.2cm}

Thus, we have proved the following.
\begin{thm}\label{thm4.1}
$(\tA(D),\td)$ is a complete metric space.
\end{thm}
\begin{prop}\label{prop4.2}
$A(D)$ is an open dense subset of $\tA(D)$. Furthermore, the
relative topology of $A(D)$ from $\tA(D)$ coincides with the usual
topology of $A(D)$.
\end{prop}
\begin{Proof}
Let $f\in A(D)$; then $f(\oD)$ is a compact subset of $\C$ and
does not contain $\infty$.

We set $\e=\frac{1}{2}\chi(f(\oD),\infty)>0$. If $g\in\tA(D)$
satisfies $\dis\sup_{|z|\le1}\chi(f(z),g(z))<\e$, then it follows
that $\infty\notin g(\oD)$. So $g\in A(D)$. This proves that
$A(D)$ is open in $\tA(D)$. Since $A(D)$ contain the polynomials
Theorem \ref{thm2.1} implies that $A(D)$ is dense in $\tA(D)$.

Let $f_n\in A(D)$ and $f\in A(D)$ be such that $f_n\ra f$
uniformly on $\oD$ with respect to the Eucledian metric in
$\C=\R^2$.  Since $\chi(\al,\bi)\le|\al-\bi|$, it follows that
$f_n\ra f$ in $\tA(D)$.

Reversely, let $g_n\in A(D)$ and $g\in A(D)$ be such that $g_n\ra
g$ in $\tA(D)$. Since $g(\oD)$ is a compact subset of $\C$, we
find $|g(z)|\le M$ for all $z\in\oD$ and some $M<+\infty$. By an
argument already used several times we have $|g_n(z)|\le M'$ for
all $z\in\oD$ and $n\ge n_5$ for some $n_5$ and $M'$, $M\le
M'<+\infty$. Thus, the convergence $g_n\ra g$ in $\tA$ implies the
convergence $g_n\ra g$ in $A(D)$. This proves that the relative
topology of $A(D)$ from $\tA(D)$ coincides with the usual topology
of $A(D)$.  \qb
\end{Proof}

In \cite{8} we consider any Hausdorff measure function $h$. We
follow the notations of \cite{8}. Let $f\in \tA(D)$. Then
$E_f=\{\zi\in T:f(\zi)\notin f(D)\}$. It is proven (\cite{8}) that
the set of $f\in A(D)$ such that $\La_h(E_f)=0$ is dense and
$G_\de$ in $A(D)$. By Proposition \ref{prop4.2} the same set is
dense and $G_\de$ in $\tA(D)$. Thus, the set of $f\in\tA(D)$ such
that $\La_h(E_f)=0$ is residual in $\tA(D)$.
\begin{prop}\label{prop4.3}
Let $h$ be any Hausdorff measure function. The set of all
$f\in\tA(D)$, such that $\La_h(E_f)=0$ is dense and $G_\de$ in
$\tA(D)$, where $E_f=\{\zi\in T:f(\zi)\notin f(D)\}$.
\end{prop}
\begin{Proof}
After the previous discussion, it remains to prove that this set
is $G_d$ in $\tA(D)$. We follow the notation of \cite{8},
especially of the first part of the proof of Theorem 3.1 of
\cite{8}. We consider the sets
$S_N=\Big\{f\in\tA(D):M_h(E_f)<\frac{1}{N}\Big\}$.

Since $\bigcap\limits^\infty_{N=1}S_N$ is our set, it suffices to
show that $S_N-\{f\equiv\infty\}$ is open. Since the singleton
$\{f\equiv\infty\}$ is a $G_\de$, the result would follow, because
the union of two $G_\de$ sets is again a $G_\de$ set.

Let $f\in S_N$ $(f\not\equiv\infty)$; thus, there exists a
countable collection of discs $D_m$ so that
$E_f\subset\bigcup\limits_mD_m$ and $\Sig\, h(r_m)<\frac{1}{N}$
where $r_m$ is the radius of $D_m$. We have $f\not\equiv\infty$.

Since $E_f\supset f^{-1}(\infty)$, by the continuity of $f$ and
the compactness of $T\sm\bigcup\limits_mD_m$, there exists
$M<+\infty$ so that $|f(\zi)|<M$ $\fa\,\zi\in T\sm\cup D_m$.

We consider a finite set of disjoint open intervals
$(\al_i,\bi_i)$ $i=1,\ld,n$,
$\al_1<\bi_1<\al_2<\bi_2<\cdots<\al_n<\bi_n<\al_1+2\pi$ and an
$r$, $0<r<1$ so that
$f^{-1}(\infty)\subset\bigcup\limits^n_{i=1}\{e^{i\vthi}:\al_i<\vthi<\bi_i\}\subset
\bigcup\limits_mD_m$ and, for every $z\in S=\Big\{te^{i\thi}:r\le
t\le1,\thi\in\bigcup\limits^n_{i=1}(\al_i,\bi_i)\Big\}$ we have
$|f(z)|>\tM$ for some $\tM$, $M<\tM<+\infty$.

Obviously, for all $\zi\in T\sm\bigcup\limits_mD_m$ we have
$f(\zi)\notin f(S)$. As $\zi\notin E_f$, it follows $f(\zi)\in
f(D\sm S)$; thus, $\de_f(f(\zi),D\sm S)>0$. By continuity of $\de_f$
and compactness of $T\sm\bigcup\limits_mD_m$ we have that
$\eta=\inf\Big\{\de_f(f(\zi),D\sm S);\zi\in
T\sm\bigcup\limits_mD_m\Big\}>0$. As in \cite{8}, if $g\in\tA(D)$
satisfies $\dis\sup_{D\sm S}|f(z)-g(z)|<\frac{\eta}{2}$ then, it
follows $\de_g(g(\zi);D\sm S)>0$ for all $\zi\in
T\sm\bigcup\limits_mD_m$. Thus $g(\zi)\in g(D\sm S)\subset g(D)$ and
$E_g\subset\bigcup\limits_mD_m$, which would imply $g\in S_N$.

So it suffices to find $\e>0$, so that
$\dis\sup_{|z|\le1}\chi(f(z),g(z))<\e$ implies $\dis\sup_{z\in D\sm
S}|f(z)-g(z)|<\frac{\eta}{2}$. Since $f(\overline{D\sm S})\subset\C$
we find $M_1<+\infty$ so that $|f(z)|<M_1$ for all
$z\in\overline{D\sm S}$. There exists $\e_1>0$ so that $|a|<M_1$,
$\chi(\al,\bi)<\e_1$ imply $|\bi|<M_1+1=M_2$. So, if
$\chi(f,g)<\e_1$, $g\in\tA(D)$ we have $|g(z)|,|f(z)|\le M_2$ for
all $z\in D\sm S$. Now for $z\in D\sm S$ we have
\[
\chi(f(z),g(z))=\frac{|f(z)-g(z)|}{\sqrt{1+|f(z)|^2}\sqrt{1+|g(z)|^2}}
>\frac{|f(z)-g(z)|}{1+M^2_2}.
\]
We set
$\e=\frac{1}{2}\min\Big(\e_1,\frac{\eta}{2(1+M^2_2)},\,\chi(f,\infty)\Big)>0$.
We can easily verify that, if $g\in\tA(D)$ is such that
$\chi(f,g)<\e$, then, $g\not\equiv\infty$ and for every $z\in D\sm
S$ we have $|f(z)-g(z)|<\frac{\eta}{2}$, which imply $g\in S_N$,
$g\not\equiv\infty$. This proves that $S_N-\{f\equiv\infty\}$ is
open in $\tA(D)$. The result easily follows since the singleton
$\{f\equiv\infty\}$ is a $G_\de$. \qb
\end{Proof}

Next we consider the sets
\[
X=\{f\in A(D):f(D)\subset f(T)\}\subset A(D)
\]
and
\[
Y=\{f\in\tA(D):f(D)\subset f(T)\}\subset\tA(D).
\]
G. Costakis pointed out that $X$ is closed in $A(D)$. In fact $Y$
is also closed in $\tA(D)$. It suffices to show that $\tA(D)\sm Y$
is open in $\tA(D)$.
\begin{prop}\label{prop4.4}
$Y$ is closed in $\tA(D)$.
\end{prop}
\begin{Proof}
Let $f\in\tA(D)\sm Y$; then, there exists $z_0\in D$ so that
$f(z_0)\in\C$ and $f(z_0)\notin f(T)$. Since $f(T)$ is a compact
subset of $\C\cup\{\infty\}$, we find $\e>0$ so that
$\chi(f(z_0),f(T))>\e$. Let $r>0$ so that $\{z:|z_0-z|\le
r\}\subset D$ and $\chi(f(z),f(z))<\e$ for all $z:|z_0-z|\le r$.
Since $f$ is non-constant we may choose $r':0<r'\le r$ so that
$f(z)\neq f(z_0)$ for all $z:0<|z-z_0|\le r'$. Let
$\de=\dis\min_{|z-z_0|=r'}|f(z)-f(z_0)|>0$. Since
$\{f(z):|z-z_0|=r'\}$ is a compact subset of $\C$, we find $\e'$,
$0<\e'<\e$, so that $\chi(f(z),w)<\e'$ for some $z:|z-z_0|=r$
implies $|f(z)-w|<\de$. We will show that if $g\in A(D)$ and
$\dis\sup_{z\in\oD}\chi(f(z),g(z))<\e'$, then $g\in\tA(D)\sm Y$.

It suffices to show that $f(z_0)\in g(D)\sm g(T)$. Since, for
every $e^{i\vthi}\in T$ we have
$\chi(g(e^{i\vthi}),f(e^{i\vthi}))<\e'$; so $\chi(g(e^{i\vthi})$,
$f(T))<\e'$; but $\chi(f(z_0),f(T))>\e>\e'$, so $g(e^{i\vthi})\neq
f(z_0)$ for all $\vthi\in\R$. It follows $f(z_0)\notin g(T)$.

We also have
\[
|f(z)-f(z_0)|\ge\de>|f(z)-g(z)| \ \ \text{for all} \ \
z:|z-z_0|=r'.
\]
Rouche's lemma implies that $f(z_0)\in g(D)$. This completes the
proof.  \qb
\end{Proof}

$\tA(D)\sm Y$ contains every non-constant polynomial. Now every
constant $a\in\C$ is the limit in $\tA(D)$ of the sequence
$P_n=a+\frac{z}{n}$. Theorem \ref{thm2.1} implies the $\tA(D)\sm
Y$ is dense; as, it is open also, it follows that $Y$ is a closed
subset of $\tA(D)$ which is of the first category. In a similar
may one can prove that $X$ is a closed subset of $A(D)$ which is
of the first category.

We notice that $X\subset Y$ are non-void. To see this, we consider
$K\subset T$ a compact set of Cantor type with zero length. It is
well known that there exists a continuous mapping $\Fi:K\ra[0,1]$
with $\Fi(K)=[0,1]$. Next we consider a Peano curve
$\Ga:[0,1]\ra\oD$ which is continuous and $\Ga([0,1])=\oD$.

We consider the function $h=\Ga\circ\Fi:K\ra\oD$, which is
continuous and $h(K)=\oD$. By the Rudin-Carleson Theorem
(\cite{7}, pages 81-82) there exists $f\in A(D)$ so that
$f_{|K}=h$ and $|f(z)|\le1$ for all $x\in\oD$. Then
$f(D)\subset\oD=h(K)=f(K)\subset f(T)$. Thus $f\in X\subset Y$ and
$X,Y$ are non-void. Summarizing we have proved
\begin{prop}\label{prop4.5}
$Y$ is a non-void closed subset of $\tA(D)$ of first category.
\end{prop}

Furthermore, we consider the set $W\subset Y$,
$W=\{f\in\tA(D):f(T)=\C\cup\{\infty\}\}\subset Y$. Obviously, $W$
is of first category in $\tA(D)$. In fact $W$ is closed in
$\tA(D)$. It suffices to show that $\tA(D)\sm W$ is open. Let
$f\in \tA(D)\sm W$. Then, there exist $P\in\C\cup\{\infty\}$ such
that $P\notin f(T)$. Since $f(T)$ is compact $\chi(P,f(T))=\de>0$.
Let $0<\e<\de$. Then, if $g\in\tA(D)$ satisfies
$\dis\sup_{z\in\oD}\chi(f(z),g(z))<\e$, it follows easily that
$P\notin g(T)$. Thus, $g\in\tA(D)\sm W$ and $\tA(D)\sm W$ is open
in $\tA(D)$. It follows that $W$ is a closed set of first category
in $\tA(D)$. We have proved
\begin{prop}\label{prop4.6}
$W$ is a closed subset of $A(D)$ of first category.
\end{prop}

We do not know if $W\neq\emptyset$. If it is true that every
compact set $K\subset T$ with zero length, is a compact of
interpolation for $\tA(D)$, then we can show that
$W\neq\emptyset$.

We consider $K_1\subset T$ a compact set of Cantor type with
length 0. $K_2\subset T$ is disjoint from $K_1$, it is a compact
set of Cantor type with length 0, which is contained in the middle
third in the first level of the construction of $K_1$. $K_2$ has
as middle point the middle point $P$ of $K$. $K_{n+1}\subset T$ is
a compact set of Cantor type with zero length, disjoint from all
$K_1,\ld,K_n$, lying in the middle third in the first level of
construction of $K_n$ and with middle point $P$.

We set $K=\{P\}\cup\bigcup\limits^\infty_{i=1}K_i$. This set is
compact with zero length. So, we assume that it is of interpolation
for $\tA(D)$. There exist continuous maps $\Fi_i:K_i\ra[2i-1,2i]$ so
that $\Fi_i(K_i)=[2i-1,2i]$. There exist Peano curves
$\Ga_i:[2i-1,2i]\ra\{z:i-1\le|z|\le i\}$, $i=1,2,\ld$ such that
$\Ga_i([2i-1,2i])=\{z:i-1\le|z|\le i\}$. We consider the maps
$\Ga_i\circ\Fi_i:K_i\ra\{z:i-1\le|z|\le i\}$ which are continuous
and onto. The map $h:K\ra\C\cup\{\infty\}$ defined by $h(P)=\infty$
and $h_{|K_i}=\Ga_i\circ\Fi_i$ is continuous and onto; thus
$h(K)=\C\cup\{\infty\}$. By assumption there exist $f\in\tA(D)$ so
that $f_{|K}=h$.

Then $f\in W$ because $f(T)\supset f(K)=h(K)=\C-\{\infty\}$. So
$W$ is non void.
\section{Further results and open questions}  
\noindent

One open question is if Mergelyan's Theorem extends to our case.
That is, to investigate the set of uniform limits with respect to
the distance $\chi$ of polynomials on a compact set $L\subset\C$
with $L^c$ connected.

The question is if the limits are exactly the continuous functions
$f:L\ra\C\cup\{\infty\}$ such that, for every compoment $V$ of
$L^0$ either $f_{|V}\equiv\infty$ or $f(V)\subset\C$ and $f_{|V}$
is holomorphic. We notice that, it is possible so that for one
component we have the first alternative, while for another
component we may have the second alternative. The general case is
open. We just give an affirmative answer in a few particular
cases.
\begin{prop}\label{prop5.1}
Let $L\subset\C$ be a compact set, such that, $L=\overline{L^0}$
and $L^0$ is starlike with respect to some point $z_0\in L^0$.
Then the uniform limits with respect to the metric $\chi$ of
polynomials on $L$ are exactly the functions $f\equiv\infty$ of
$f:L\ra\C\cup\{\infty\}$ continuous with $f(L^0)\subset\C$ and
$f_{|L^0}$ holomorphic.
\end{prop}

For the proof we may assume $z_0=0$. Then we imitate the proof on
$L=\oD$. The difference is that when we approximate $f$ on the
compact set $rL$, with $0<r<1$, we cannot use a Taylor development
of $f$, but we can use the classical Mergelyan's Theorem. The
approximation is uniform on $rL$ with respect to the Eucledian
distance on $\C=\R^2$; since $\chi(\al,\bi)\le|\al-\bi|$ for all
$\al,\bi\in\C$, this implies approximation uniform with respect to
the distance $\chi$.

Under the assumptions of Proposition \ref{prop5.1}, if $w\in\C\sm
L$, then the uniform limits on $L$ by functions of the form
$P\Big(\frac{1}{z-w}\Big)$, where $P$ is any polynomial are the
same with those given by Proposition \ref{prop5.1}. The difference
in the proof is, when we apply Mergelyan's Theorem to do
approximation of $f$ on $rL$, $0<r<1$, then we choose as pole
$rw\notin rL$ and not $\infty$.

Next we consider the analogue problem on an annulus. Without loss
of generality we may assume that this annulus is
$\OO=\cD(0,r,1)=\{z\in\C:r<|z|<1\}$ with $0<r<1$.
\begin{prop}\label{prop5.2}
Let $0<r<1$ and $\OO=\cD(0,r,1)$. Then the uniform limits on $\oO$
with respect to $\chi$ of polynomials are the functions
$f\in\tA(D)$.
\end{prop}

In the proof, if the sequence of polynomials $P_n$ converges
uniformly on $\overline{\cD(0,r,1)}$ with respect to $\chi$ to
$f\not\equiv\infty$, then $f\{z:r<|z|<1\}\subset\C$ and $f$ is
holomorphic in $\{z:r<|z|<1\}$. So there exists $M<+\infty$, so
that $|f(z)|<M$ on $|z|=\frac{r+1}{2}$. Because $\chi(P_n,f)\ra0$
we find $M_1:M<M_1<+\infty$ and $n_0$ so that $|P_n(z)<M_1$ for
$|z|=\frac{1+r}{2}$ and $n\ge n_0$. For $|a|,|b|<M_1$ we have
$\chi(\al,\bi)\ge\frac{|al-\bi|}{1+M^2_1}$; So $P_n\ra f$
uniformly on $|z|=\frac{1+r}{2}$ with respect to the usual
Euclidian metric on $\C=\R^2$. Thus, the sequence $P_n$,
$n=1,2,\ld$ is uniformly Cauchy on $|z|=\frac{r+1}{2}$ with
respect to the Eucledian metric in $\C$. The maximum principle
implies that the sequence of polynomials $P_n$ is uniformly Cauchy
on $|z|\le\frac{r+1}{2}$. Its limit is a holomorphic extension of
$f$, so finally $f\in\tA(D)$. It is remarkable that the value
$\infty$ is not permitted on $|z|=r$, unless $f\equiv\infty$.

If we wish to investigate the analogue problem on
$\overline{\cD(0,r,1)}$, $0<r<1$ using functions of the form
$P\Big(\frac{1}{z-w}\Big)$ where $P$ is any polynomial and
$w\in\C$ fixed, then the answer is the following.

If $|w|>1$ then the possible limits are $f\in\tA(D)$. The proof is
the same as in the case of polynomials.

The case $|w|<r$ is reduced to the previous one by an inversion. The
possible limits are $f\equiv\infty$ or $f:\{z:r\le|z|\}\ra
C\cup\{\infty\}$ continuous with $f(\{z:r<|z|\})\subset\C$ and $f$
holomorphic in $r<|z|$ and $\dis\lim_{z\ra\infty}f(z)$ exist in
$\C$.

In the case $r<|w|<1$, the limits are exactly $f\equiv\infty$ or
$f(z)=P\big(\frac{1}{z-w}\Big)$ with $P$ any polynomial. The
argument in the proof uses Cauchy transforms as in \S\,2. In the
case $|w|=1$ or $|w|=r$ I do not know the answer.

For $r=1$, I do not know the answer to the analogue question of
finding the uniform limits on $T=\{z\in\C:|z|=1\}$ with respect to
$\chi$ of all polynomials.

Finally we may consider $\oO=\overline{\cD(0,r,1)}$ $0<r<1$ and
searching for the limits using two poles. There are 5 cases
$\{\infty\}\cup\{w\in\C:|w|>1\}$, $\{w\in\C:|w|=1\}$,
$\{w\in\C:r<w<1\}$, $\{w\in\C:|w|=r\}$, $\{w\in\C:|2|<r\}$. So the
number of possible locations of two poles $w_1,w_2$ is
$\binom{5}{2}=10$. In some of these cases the author does not know
the answer, but in some other cases we know the answer.

If we use the poles $w_1=\infty$, $w_2=0$ then we approximate by
functions of the form $\ssum^N_{n=-N}a_nz^n$. The limit functions
are $f\equiv\infty$ or
$f:\overline{\cD(0,r,1)}\ra\C\cup\{\infty\}$ continuous,
$f(\cD(0,r,1))\subset\C$ and $f_{|\cD(0,r,1)}$ holomorphic.

For the proof we use Laurent expansion and we approximate in
$\chi$ metric separately the principal part and the regular part.
Next we need a lemma of the following form.
\begin{lem}\label{lem5.3}
Let $g,f:\cD(0,r,1)\ra\C\cup\{\infty\}$ are continuous on
$\overline{\cD(0,r,1)}$, $f(\{z:r\le|z|<1\})\subset\C$,
$g(\{z:r<|z|\le1\})\subset\C$,
$f_n,g_n:\overline{\cD(0,r,1)}\ra\C$ continuous.

We assume that
\[
\sup_{z\in\overline{\cD(0,r,1)}}\chi(f_n(z),f(z))\ra0 \ \
\text{and} \ \
\sup_{z\in\overline{\cD(0,r,1)}}\chi(g_n(z),g(z))\ra0.
\]
Then $f+g$ is well defined with values in $\C\cup\{\infty\}$ and
\[
\sup_{z\in\overline{\cD(0,r,1)}}\chi(f_n(z)+g_n(z),f(z)+g(z))\ra0.
\]
\end{lem}

For the proof of the lemma it is essential that
$f^{-1}(\infty)\subset\{z:|z|=1\}$ and
$g^{-1}(\infty)\subset\{z:|z|=r\}$; so the compact sets
$f^{-1}(\infty)$, $g^{-1}(\infty)$ are disjoint.

We fix $\e>0$ and we are looking for $n_0$ so that
$\dis\sup_{z\in\overline{\cD(0,r,1)}}\chi(f_n(z)+g_n(z),f(z)+g(z))<\e$
for all $n\ge n_0$.

We obtain that there exist a constant $M<+\infty$ so that
$|f(z)|\le M$ for $r\le|z|\le\frac{r+1}{2}$ and $|g(z)|\le M$ for
$\frac{r+1}{2}\le|z|\le1$. We can choose a relatively open
neighborhood $V_1$ of $f^{-1}(\infty)$ so that
$V_1\subset\{z:\frac{r+1}{2}<|z|\le1\}$ so that $|f(z)|>M_1$ for
$z\in V_1$, where $M_1$ is in our disposal. Then
$|f(z)+g(z)|>M_1-M$ on $V_1$.

We easily find $n_1$ so that for $n\ge n_1$ we have
$|f_n(z)+g_n(z)|>M_1-M-1$ on $V_1$. We can choose $M_1$ big enough
so that the $\chi$ diameter of
$\{\infty\}\cup\{w\in\C:|w|>M_1-M-1\}$ is less than $\e$. So, for
$n\ge n_1$ we have $\dis\sup_{z\in
V_1}\chi(f_n(z)+g_n(z),f(z)+g(z))<\frac{\e}{3}$. In a similar way we
find a relatively open neighborhood $V_2$ of $g^{-1}(\infty)$,
$V_2\subset\Big\{z:r\le|z|<\frac{r+1}{2}\Big\}$ and $n_2$ so that
for $n\ge n_2$ we have $\dis\sup_{z\in
V_2}\chi(f_n(z)+g_n(z),f(z)+g(z))<\frac{\e}{3}$. The set
$\overline{\cD(0,r,1)}-(V_1\cup V_2)$ is compact. By continuity
there exist $M_2<+\infty$ so that $|f(z)|,|g(z)|<M_2$ on
$\overline{\cD(0,r,1)}\sm(V_1\cup V_2)$. One easily finds $n_3$ so
that $|f_n(z)|,|g_n(z)|<M_2+1$ on $\overline{\cD(0,r,1)}\sm(V_1\cup
V_2)$ for $n\ge n_3$. We set $M_3=M_2+1$. Since for
$|\al|,|\bi|<M_3$ we have $\chi(\al,b)\ge\frac{|\al-\bi|}{1+M^2_3}$,
we see that $f_n\ra f$ and $g_n\ra g$ uniformly on
$\overline{\cD(0,r,1)}\sm(V_1\cup V_2)$ with respect to the
Eucledian metric in $\C=\R^2$. Thus, $f_n+g_n\ra f+g$ uniformly on
$\overline{\cD(0,r,1)}\sm(V_1\cup V_2)$ with respect to the
Eucledian metric in $\C=\R^2$. Since $\chi(\al,\bi)\le|\al-\bi|$, it
follows that $f_n+g_n\ra f+g$ uniformly on
$\overline{\cD(0,r,1)}\sm(V_1\cup V_2)$ with respect to $\chi$
metric. Thus, there exists $n_4$, so that
$\dis\sup_{z\in\overline{\cD(0,r,1)}\sm(V_1\cup
V_2)}\chi(f_n(z)+g_n(z),f(z)+g(z))<\frac{\e}{3}$. We set
$n_0=\max(n_1,n_2,n_4)$ and we have the result. \qb\vspace*{0.2cm}

The limiting case $r=1$ of the previous result is quite different.
\begin{thm}\label{thm5.4}
The uniform limits of the trigonometric polynomials
$\ssum^N_{-N}a_nz^n$ on $T=\{z:|z|=1\}$ with respect to the $\chi$
metric are exactly all continuous functions
$f:T\ra\C\cup\{\infty\}$.
\end{thm}
\begin{Proof}
One direction is obvious. Let $f:T\ra C\cup\{\infty\}$ be a
continuous function and $\e>0$. We have to find a trigonometric
polynomial $Q$, $Q(z)=\ssum^N_{-N}a_nz^n$, so that
$\dis\sup_{|z|=1}\chi(f(z),Q(z))<\e$. If $f\equiv\infty$ it suffices
to choose $Q(z)=nz$ for $n$ big enough.

Assume there exists $z_0\in T$ so that $f(z_0)\in\C$. If
$f(z)\neq\infty$ for all $z\in T$, then $f:T\ra\C$ being
continuous it can be approximated by a trigonometric polynomial
$Q$ so that $\dis\sup_{z\in T}|f(z)-Q(z)|<\e$. Since
$\chi(\al,\bi)\le|\al-\bi|$, it follows $\dis\sup_{z\in
T}\chi(f(z),Q(z))<\e$.

It remains to examine the case $f(z_0)\in\C$ and $f(z_1)=\infty$ for
some $z_0,z_1\in T$. We consider $M<+\infty$ so that the $\chi$
diameter of the set $E=\{\infty\}\cup\{z\in\C:|z|\ge M\}$ is less
than $\frac{\e}{2}$. Then $f^{-1}(E^0)$ is an open subset of $T$
containing $f^{-1}(\infty)$. The compact set $f^{-1}(\infty)$ is
covered by a finite number of components of $f^{-1}(E)$.

Thus,
$f^{-1}(\infty)\subset\Big\{e^{i\thi}:\thi\in\bigcup\limits^L_{n=1}(\al_n,\bi_n)\Big\}=V$
where $\al_1<\bi_1<\al_2<\bi_2<\cdots<\al_L<\bi_L<\al_1+2\pi$ and
$|f(e^{i\al_n})|=|f(e^{9\bi_n})|=M$ for all $n=1,\ld,L$. We define
a new function $g:T\ra\C$ setting $g(e^{i\vthi})=f(e^{i\vthi})$
for $e^{i\vthi}\in T\sm V$. For each $n=1,\ld,L$ we choose
$g_{|\{e^{i\thi}:\thi\in(\al_n,\bi_n)\}}:\{e^{i\thi}:\thi\in(\al_n,\bi_n)\}\ra
\{z\in\C:|z|=M\}$ continuous so that
$\dis\lim_{\thi\ra\al^+_n}g(e^{i\thi})=f(e^{i\al_n})$ and
$\dis\lim_{\thi\ra\bi^-_n}g(e^{i\thi})=f(e^{i\bi_n})$. For this,
it suffices to follow one of the arcs of the circle
$\{z\in\C:|z|=M\}$ with starting point $f(e^{i\al_n})$ and ending
to $f(e^{i\bi_n})$. It is easy to check that $g:T\ra\C$ is
continuous and that, by the choice of $M$ we have $\dis\sup_{z\in
T}\chi(f(z),g(z))<\frac{\e}{2}$. Now $g$ taking values on $\C$ and
being continuous it can be approximated uniformly on $T$ by a
trigonometric polynomial $Q$, with respect to the Eucledian metric
on $\C=\R^2$. We choose $Q(z)=\ssum^N_{-N}a_nz^n$ so that
$\dis\sup_{z\in T}|g(z)-Q(z)|<\frac{\e}{2}$. Since
$\chi(\al,\bi)\le|\al-\bi|$ it follows $\dis\sup_{z\in
T}\chi(g(z),Q(z))<\frac{\e}{2}$. Since we also have
$\dis\sup_{z\in T}\chi(f(z),g(z))<\frac{\e}{2}$ the triangular
inequality implies $\dis\sup_{z\in T}(f(z),Q(z))<\e$ and the proof
is completed.  \qb
\end{Proof}

If $I$ is a compact segment or a homeomorphic image of a compact
segment, $I\subset\C$, then the uniform on $I$ limits of the
polynomials are all continuous functions $f:T\ra\C\cup\{\infty\}$.
This is another particular case of the required extension of
Mergelyan's Theorem. The proof is similar to that of Theorem
\ref{thm5.4} with the difference that we do not approximate by a
trigonometric polynomial but by a polynomial. This is possible by
the classical Mergelyan's Theorem, because $I^0=\emptyset$ and $I^c$
is connected.

Another particular case of possible extension of Mergelyan's
Theorem is the following.

Let $L_1,L_2,\ld,L_N$ be a finite collection of compact subsets of
$\C$ with connected complements. We assume that for each $L_i$,
$i=1,\ld,N$ the extension of Mergelyan's Theorem is valid. We also
assume that the $L_i$'s, $i=1,\ld,N$ are two by two disjoint. Then
the extension of Mergelyan's Theorem is valid for $L_1\cup
L_2\cup\cdots\cup L_N$. Indeed, if $f:L_1\cup\cdots\cup
L_N\ra\C\cup\{\infty\}$ is continuous and for every $i=1,\ld,N$
either $f_{|L_i}\equiv\infty$ or $f(L^0_i)\subset\C$ and
$f_{|L^0_D}$ holomorphic, then by assumption there exist
polynomials $P_i$, $i=1,\ld,N$ so that $\dis\sup_{z\in
L_i}\chi(f(z),P_i(z))<\frac{\e}{2}$ for all $i=1,\ld,N$. Now by
the classical Mergelyan's Theorem on $L=L_1\cup L_2\cup\cdots\cup
L_N$, there exists a polynomial $P$ so that $\dis\sup_{z\in
L_i}|P(z)-P_i(z)|<\frac{\e}{2}$. Since
$\chi(\al,\bi)\le|\al-\bi|$, we have $\dis\sup_{z\in
L_i}\chi(P(z),P(z)|<\frac{\e}{2}$.

The triangular inequality gives the result.

So in particular the extension of Mergelyan's Theorem is valid for
the union of two disjoint closed discs. What happens for the union
of a sequence of pairwise disjoint closed discs converging to one
point $P$ in $\C$? The difficulty is when $f(P)=\infty$.

\newpage
\begin{center}
{\bf {\Large Appendix}}
\end{center}

In this Appendix I collect the basic facts about the metric $\chi$
we used (and proved) in this article.

Because $\chi(\al,\bi)\le|\al-\bi|$ for all $\al,\bi\in\C$ it
follows that uniform convergence with respect to the Eucledian
metric in $\C=\R^2$ implies uniform convergence with respect to
the metric $\chi$. The converse does not hold in general. However,
if $f_n\ra f$ uniformly with respect to $\chi$ on some set $E$,
and there exist $M<+\infty$ so that $|f_n(z)|\le M$ and $|f(z)|\le
M$ for all $z\in E$ an all $n\ge n_0$ for some $n_0$, then it
follows that $f_n\ra f$ uniformly on $E$ with respect to the
Eucledian metric on $\C=\R^2$. The reason is that
$\chi(f_n(z),f(z))\ge\frac{|f_n(z)-f(z)|}{1+M^2}$.

In fact it suffices that only $|f(z)|\le M$ on $E$ (or that only
$|f_n(z)\le M$ on $E$ for all $n\ge n_0$ for some $n_0$). The
reason is that when $\dis\sup_{z\in E}\chi(f_n(z),f(z))\ra0$ and
$|f(z)|\le M$ on $E$, then there exists $n_0$ so that for $n\ge
n_0$ we have $|f_n(z)|\le M+1$ on $E$.

Another fact we have used (and proved) is that, if $\OO$ is a
domain and $f_n:\OO\ra\C$ a sequence of holomorphic functions in
$\OO$ and there exists a function $f:\OO\ra\C\cup\{\infty\}$ so
that $\dis\sup_{z\in\OO}\chi(f_n(z),f(z))\ra0$, then either
$f\equiv\infty$ or $f(\OO)\subset\C$ and $f$ is holomorphic in
$\OO$. This follows by Hurwitz Theorem.

Finally we have proven the following on an annulus
$\overline{D(0,r,1)}$, $0<r<1$.

If $f_n:\overline{\cD(0,r,1)}\ra\C$ and
$g_n:\overline{\cD(0,r,1)}\ra\C$ and
$f:\overline{\cD(0,r,1)}\ra\C\cup\{\infty\}$ and
$g:\overline{\cD(0,r,1)}\ra C\cup\{\infty\}$ are given. $f$ and
$g$ are continuous and $f^{-1}(\infty)\subset\{z:|z|=1\}$ and
$g^{-1}(\infty)\subset\{z:|z|=r\}$.

Assume $\dis\sup_{z\in\overline{\cD(0,r,1)}}\chi(f_n(z),f(z))\ra0$
and $\dis\sup_{z\in\overline{\cD(0,r,1)}}\chi(g_n(z),g(z))\ra0$.

Then
$\dis\sup_{z\in\overline{\cD(0,r,1)}}\chi(f_n(z)+g_n(z),f(z)+g(z))\ra0$
and $f+g$ is well defined.

I think the above facts except probably the last one are well known.
May be we can find them in the works of Caratheodory, or Zalcman.
The proofs may be shortened if we did references to
articles proving the above facts.\vspace*{0.3cm} \\
{\bf Acknowledgment.} I would like to thank G. Costakis for
bringing to my attention reference \cite{8} and for a helpful
communication.
\vspace*{1cm}
University of Athens \\
Department of Mathematics  \\
157 84 Panepistemiopolis \\
Athens \\
GREECE\\
e-mail: vnestor@math.uoa.gr

\end{document}